\documentclass {article}

\usepackage{amsfonts,amsmath,latexsym}
\usepackage{amstext,amssymb}
\usepackage{amsthm}
\usepackage{verbatim}

\newcommand {\R}{\mathbb{R}}
\newcommand {\Ss}{\mathbb{S}}

\newcommand {\del}{\partial}
\newcommand {\Div}{\operatorname{div}}

\newcommand {\vol}{\operatorname{Vol}}

\newcommand {\calS}{\mathcal{S}}

\newcommand {\calW}{\mathcal{W}}

\newtheorem {thm} {Theorem}
\newtheorem {prop} [thm] {Proposition}
\newtheorem {lemma} [thm] {Lemma}

\newtheorem {rmk} {Remark}

\begin {document}

\title{ 
{\bf Eigenvalues of Euclidean wedge 
domains in higher dimensions}
}
\author{Jesse Ratzkin} 
\maketitle

\begin {abstract} 
In this paper, we use a weighted isoperimetric 
inequality to give a lower bound for the first
Dirichlet eigenvalue of the Laplacian on a 
bounded domain inside a Euclidean cone. Our 
bound is sharp, in that only sectors realize it. 
This result generalizes a lower bound of Payne and 
Weinberger in two dimensions. 
\end {abstract}

\section{Introduction}

Lower bounds for the first Dirichlet 
eigenvalue of the Laplacian often arise 
from an integral 
inequality relating the boundary of a
domain to its interior (see, for example,
 \cite{Ch, LT, P}). 
Moreover, these inequalities can lead one 
to a characterization of the optimal
domains, that is, the domains 
for which the inequality is an equality. 
The Faber-Krahn inequality \cite{F,K}
provides the model for such a lower bound, 
using the classical isoperimetric inequality 
to prove that, among all domains with the 
same volume in a simply connected space of 
constant curvature, geodesic balls have the
least eigenvalue. In 
this paper we prove a similar estimate
for domains inside a convex cone. 

To state our main theorem, we first introduce some 
notation. Let $\Omega$ be a convex domain in the upper 
unit hemisphere $\Ss^{n-1}_+ = \Ss^{n-1} \cap 
\{x_n > 0\}$, and let 
$$\calW = \{(r,\theta) \mbox{ }|\mbox{ } r\geq 0, 
\theta \in \Omega \}$$
be the cone over $\Omega$. For $r_0 > 0$, define the 
sector 
$$\calS_{r_0} = \{(r,\theta) \mbox{ }|\mbox{ } 0 \leq r 
\leq r_0, \theta \in \Omega\}.$$
If $D$ is a bounded domain in $\calW$, we denote 
the first Dirichlet eigenvalue of the Laplacian on 
$D$ by $\lambda_1(D)$. 

Next, let $\mu$ be the first Dirichlet  eigenvalue of 
the Laplacian on 
$\Omega$, with eigenfunction $\psi$. Normalize
$\psi$ so that $\int_{\Omega} \psi^2 = 1$. 
Observe that
\begin {equation} \label{defn-harm-weight}
 w(r,\theta) = r^\alpha \psi(\theta), \qquad
\alpha := \frac{2-n}{2} + \sqrt{ \left ( 
\frac{2-n}{2} \right )^2 + \mu} \end {equation}
is a positive harmonic function in $\calW$, which 
is zero on the boundary of $\calW$. Monotonicity, 
with $\Omega \subset \Ss^{n-1}_+$, implies
$\mu > n-1$, and so $\alpha > 1$. 

\begin {thm} \label{main-thm}
Let $D$ be a bounded domain in the cone $\calW$, 
and choose $r_0$ so that 
$$\int_D w^2 dV = \int_{\calS_{r_0}} w^2 dV.$$
Then $\lambda_1(D) \geq \lambda_1(\calS_{r_0})$, 
with equality if and only if $D = \calS_{r_0}$ 
almost everywhere. 
\end {thm}

This theorem extends an eigenvalue estimate 
of Payne and Weinberger \cite{PW} to dimensions 
greater than two. 
Like their estimate, ours relies on 
a weighted isoperimetric inequality; such 
inequalities are sometimes called isoperimetric 
inequalities for measures \cite{Ro} or 
isoperimetric inequalities for densities 
\cite{M, RCBM}. 

The rest of this paper proceeds 
as follows. First we compute the first Dirichlet 
eigenvalue of the Laplacian on a sector in 
Section \ref{sector-section}, 
both in terms of the radius $r_0$ and the integral of 
the weight function $w$. The bulk of this paper 
lies in Section \ref{isop-sec}, where we prove 
a weighted isoperimetric inequality for domains 
in the cone $\calW$. Finally, in Section 
\ref{rayleigh-sec} we use our isoperimetric 
inequality to estimate the Rayleigh quotient, 
proving our main theorem. 

We conclude this introduction by briefly 
comparing the two-dimensional and 
higher-dimensional theorems. First observe
that in the two-dimensional 
case, the link $\Omega$ is an interval 
$(0,\pi/\alpha)$, for some $\alpha>1$, and the 
cone $\calW$ has the form 
$$\calW = \{ (r, \theta) \mbox{ }|\mbox{ } 0 
< \theta < \pi/\alpha \}$$
in polar coordinates. 
In both the two-dimensional and 
higher-dimensional cases, 
the key weighted isoperimetric inequality 
comes from an inequality for domains in a 
half-space (see Lemma \ref{half-space-isop}) 
below). In two dimensions, there is an obvious 
way to open the cone up to a half-plane, 
whereas in higher dimensions this is more 
subtle. The second  main 
difference between the two proofs is 
a technical complication. In the two-dimensional 
case, the first eigenvalue of the link 
is $\mu = \alpha^2$, and the harmonic weight is 
$w(r,\theta) = r^\alpha \sin (\alpha\theta)$ (up to 
a constant multiple). The fact the the opening 
angle of the cone, the first eigenvalue of the 
link, and the exponent of the radial part of the 
harmonic weight function all agree simplifies 
much of the analysis. 

\medskip
\small{\noindent  {\bf Acknowledgements.} This 
work is 
supported by the Science Foundation of Ireland 
under grant number MATF636. I would also like to 
thank N. Korevaar, M. Ritor\'e, A. Ros, and A. 
Treibergs for helpful conversations. }

\section{First eigenvalue of a sector}
\label{sector-section}

In this section we compute the first Dirichlet 
eigenvalue of a sector $\calS_{r_0}$, to compare
with our estimates for the Rayleigh quotient 
below in Section \ref{rayleigh-sec}. 

\begin {lemma} \label{eigen-sector}
Let 
$$a = \sqrt{\left ( \frac{n-2}{2} \right )^2
+ \mu} = \alpha + \frac{n-2}{2}$$
and consider the sector $\calS_{r_0}$ for 
$r_0> 0$. Then 
$$\lambda_1 ( \calS_{r_0}) = 
\frac{j_a^2}{r_0^2} = \left [ 
(2a+2) \int_{\calS_{r_0}} w^2 dV
\right ]^{-\frac{1}{a+1}} j_a^2,$$
where $j_a$ is the first positive 
zero of the Bessel function $J_a$. 
\end {lemma}

Notice that $a=\alpha$ in dimension two. 

\begin {proof} 
We write the eigenfunction for the sector 
as $u(r,\theta) = f(r) \psi(\theta)$, and denote 
eigenvalue by $\lambda$. Then we rewrite the 
eigenvalue equation  
$$-\lambda u = \Delta u = u_{rr} + \frac{n-1}{r} 
u_r + \frac{1}{r^2} \Delta_{\Ss^{n-1}}u,$$
as 
\begin{equation} \label{sector-ode1}
0 = r^2 f'' + (n-1) rf' + (\lambda r^2 
- \mu) f.
\end {equation} 
Look for a solution of the form $f(r) 
= r^q J_a(kr)$, where $J_a$ is the Bessel 
function of weight $a$, and $k$ and $q$ 
are positive constants. Observe that 
$$k^2r^2 J_a''(kr) + krJ_a'(kr) + 
(k^2 r^2 - a^2) J_a(kr) = 0.$$
Plugging this choice of $f$ into 
equation (\ref{sector-ode1}), we get
\begin{eqnarray*}
0 & = & r^2 f'' + (n-1) rf' + (\lambda r^2 
- \mu) f\\
& = & q(q-1) r^q J_a + 2kqr^{q+1}
J_a' + k^2r^{q+2} J_a'' + 
(n-1)qr^q J_a + (n-1) kr^{q+1}J_a' 
+ (\lambda r^2- \mu)r^q J_a(kr) \\
& = & r^q [k^2 r^2 J_a'' + k(2q+n-1) rJ_a'
+(\lambda r^2-\mu +q(n+q-2))J_a]\\
& = & r^q [k(2q+n-2)rJ_a' 
+ ((\lambda - k^2)r^2 - \mu + q(n+q-2) 
+a^2)J_a].
\end{eqnarray*}
Setting the individual terms to zero, we
obtain 
\begin {equation} \label{sector-param}
q = \frac{2-n}{2}, \qquad k= \sqrt{\lambda}, 
\qquad a^2 = \mu -q(n+q-2) = \mu +\left (
\frac{n-2}{2} \right )^2.
\end {equation}
The fact that the eigenfunction 
vanishes on the boundary of $\calS_{r_0}$ 
forces 
$$0 = f(r_0) = r_0^{\frac{2-n}{2}} J_a
(\sqrt{\lambda} r_0) \Rightarrow 
\sqrt \lambda r_0 = j_a.$$
In other words, the eigenvalue is 
$$\lambda = \lambda_1(\calS_{r_0}) 
= \frac{j_a^2}{r_0^2}, \qquad 
a = \sqrt{\left ( \frac{n-2}{2} \right )^2
+ \mu}$$
and the eigenfunction is 
$$u(r,\theta) = r^{\frac{2-n}{2}} J_a
\left ( \frac{j_a r}{r_0} \right ) 
\psi(\theta).$$

Finally, we rewrite our expression for 
$\lambda$ in terms of the integral 
$\int_{\calS_{r_0}} w^2 dV$. Recall
$$\alpha = \frac{2-n}{2} + \sqrt{ 
\left( \frac{2-n}{2} \right )^2 + \mu}
= a-\frac{n-2}{2}.$$
By explicit computation, 
$$\int_{\calS_{r_0}} w^2 dV = 
\int_\Omega \psi^2 dA_{\Ss^{n-1}}
\int_0^{r_0} r^{2\alpha+n-1} dr = 
\frac{r_0^{2\alpha + n}}{2\alpha + n}
= \frac{r_0^{2a+2}}{2a+2},$$
or 
$$r_0 = \left [ (2\alpha + n) 
\int_{\calS_{r_0}} w^2 dV \right ]
^{\frac{1}{2\alpha + n}} = 
\left [ (2a+2) 
\int_{\calS_{r_0}} w^2 dV \right 
]^{\frac{1}{2a+2}}.$$
Plugging this into our formula for 
$\lambda$, we see
$$\lambda_1 (\calS_{r_0}) = \left [
(2a+2) \int_{\calS_{r_0}} w^2 dV \right 
]^{-\frac{1}{a+1}} j_a^2.$$
\end {proof}

\section{Isoperimetric inequality}
\label{isop-sec}

The goal of this section is to establish 
the following weighted isoperimetric 
inequality for bounded domains $D \subset 
\calW$. 
\begin {prop} 
Let $D \subset \calW$ be a bounded domain 
with piecewise smooth boundary.  
Then 
\begin {equation} \label {weight-isop}
\int_{\del D} w^2 dA \geq \left [ (2a+2)
\int_D w^2 dV \right ]^{\frac{2a+1}{2a+2}}, 
\qquad a = \sqrt{\left ( \frac{n-2}{2} 
\right )^2 + \mu} = \alpha + \frac{n-2}{2}.
\end {equation}
Moreover, equality implies $D$ is (almost 
everywhere) a sector $\calS_{r_0}$ for some 
$r_0 > 0$. 
\end {prop}
\begin {rmk} In proving our weighted isoperimetric 
inequality, we will not explicitly compute the 
coefficient $(2a+2)^{\frac{2a+1}{2a+2}}$, but 
rather show that there is some constant $c$ 
depending only on $n$ and $\Omega$ such that 
$$\int_{\del D} w^2 dA \geq c \left [ \int_D
w^2 dV \right ]^{\frac{2a+1}{2a+2}}.$$
Then we recover the constant $c$ by computing 
the relevant integrals for a sector. 
\end {rmk}

We start with the following weighted 
isoperimetric problem in a half-space.  

\begin{lemma} \label{half-space-isop}
Let $\widetilde D 
\subset \R^n_+ = \{ x_n \geq 0\}$ be a 
bounded domain with piecewise smooth 
boundary. If we hold $\int_{\widetilde 
D} x_n^2 dV$ constant, the minimizers of 
$\int_{\del \widetilde D} x_n^2 dA$ are 
hemispheres whose centers lie on the 
hyperplane $\{x_n = 0\}$. 
\end {lemma}

\begin {proof} If $n=1$ the lemma is 
automatically true. Payne and 
Weinberger \cite{PW} proved that in 
dimension two the minimizers of the variational
problem are semi-circles centered on the 
$x_1$ axis. Indeed, if one lets $\widetilde D = 
\{ (x_1, x_2) \mbox{ }|\mbox{ } 0 \leq x_2 
\leq \phi(x_1) \}$, a straight-forward 
computation shows that the 
Euler-Lagrange equation is 
\begin{equation} \label{euler-lagrange-ode}
2 \sqrt{1+(\phi')^2} -
\frac{2(\phi')^2}{\sqrt{1+(\phi')^2}}
- \phi\left ( \frac{\phi'}{\sqrt {1+(\phi')^2}} 
\right )' = \Lambda \phi, \end {equation}
where the constant $\Lambda$ is the Lagrange 
multiplier. For any choice of constants $R$ 
and $c$, the function $\phi(x_1) = 
\sqrt{R^2 - (x_1 - c)^2}$ solves equation 
(\ref{euler-lagrange-ode}), so, by uniqueness 
of solutions to ODEs, semi-circles are the 
only critical points and thus minimize. We 
pause here to remark that the Euler-Lagrange 
equation in higher dimensions is 
$$ 2\sqrt{1+|\nabla \phi |^2} - \frac{2 |\nabla 
\phi |^2}{\sqrt{1+|\nabla \phi |^2}} - \phi 
\Div \left ( \frac{\nabla \phi}{\sqrt{1+ |\nabla 
\phi}} \right ) = \Lambda \phi.$$

For $n\geq 3$ we reduce the variational 
problem to a two-dimensional problem using 
Steiner symmetrization (see note A of 
\cite{PS}), starting with $\{ x_1 = 0 \}$ 
as the symmetry hyperplane. Given a point 
$\bar x = (0, x_2, \dots, x_n)$
in the symmetry hyperplane, let $l_{\bar x}$ 
be the line through $\bar x$ perpendicular to 
to $\{ x_1 = 0\}$. By Sard's theorem, for 
almost every $\bar x$ such that the intersection
$\widetilde D \cap l_{\bar x}$ is nonempty one 
can write 
\begin {equation} \label{graphs-1}
\widetilde D \cap l_{\bar x} = 
\cup_{k=1}^{m_{\bar x}}
((a_k, x_2, \dots, x_n), 
(b_k, x_2, \dots, x_n)).
\end{equation} 
The symmetrized domain $\widetilde D_1$ is 
given by 
\begin {equation} \label{graphs-2}
\widetilde D_1 \cap l_{\bar x} = ((-\bar a, 
x_2, \dots, x_n),(\bar a, x_2, 
\dots, x_n)), \qquad 2\bar a = 
\sum_{k=1}^{m_{\bar x}} b_k - a_k.
\end {equation}
Observe that $\widetilde D_1$ is now symmetric 
about the $\{ x_1 = 0\}$ hyperplane, and has 
the same volume as $\widetilde D$. Also, because 
$x_n$ is constant along the line $l_{\bar x}$ 
and the length of $\widetilde D \cap l_{\bar x}$
is the same as the length of $\widetilde D_1 
\cap l_{\bar x}$, Fubini's theorem implies 
$$\int_{\widetilde D_1} x_n^2 dV = \int_{\widetilde D}
x_n^2 dV.$$

Next we show examine the boundary integral. 
Let $G$ be the projection of $\del 
\widetilde D$ onto the $\{ x_1 = 0\}$ 
hyperplane, and notice that $G$ is also 
the orthogonal projection of $\widetilde D_1$ onto 
$\{ x_1 = 0\}$. Moreover, our representations 
of $\widetilde D$ and $\widetilde D_1$ as 
graphs over the $\{ x_1 = 0 \}$ hyperplane 
in equations (\ref{graphs-1}) and 
(\ref{graphs-2}) are smooth almost everywhere, 
and so the integral identities below are valid. 
Now rewrite the boundary integral as 
\begin {eqnarray*}
\int_{\del \widetilde D} x_n^2 dA 
& = & \int_0^\infty \int_{\del \widetilde D \cap 
\{x_n = t\} } x_n^2 dA \\
& = & \int_0^\infty t^2 \int_{G \cap 
\{x_n = t\}}  \sqrt {1+ \sum_{k=1}
^m\sum_{j=2}^n \left (\frac
{\del a_k}{\del x_j} \right )^2 + 
\left ( \frac{\del b_k}{\del x_j}
\right )^2} dx_2\cdots dx_{n-1} dt \\
& \geq & \int_0^\infty t^2 
\int_{G \cap \{ x_n = t\}} 
\sqrt{1+ 2 \sum_{j=2}^n \left (
\frac{\del \bar a}{\del x_j} \right )^2}
dx_2\cdots dx_{n-1} dt \\
& = & \int_{\del \widetilde D_1} x_n^2 dA. 
\end {eqnarray*}
Here we have used the fact that $2\bar a 
= \sum_{k=1}^m b_k - a_k$ and Minkowski's 
inequality. 

The lemma now follows by induction. 

\end {proof}

\begin {rmk} 
One can adapt this proof to show that for any 
$q>-1$ the minimizers of $\int_{\del \widetilde 
D} x^q dA$, with $\int_{\widetilde D} x_n^q dV$
held constant, are also hemispheres whose 
centers lie on the $\{ x_n = 0 \}$ hyperplane. (The 
Euler-Lagrange equation is almost the same for the 
case of an arbitrary power $q$.) The 
case of $q=0$ is the classical isoperimetric 
inequality. The condition $q>-1$ ensures that 
the relevant integrals converge. 
\end {rmk}

Let 
$$\omega_k = \vol(\Ss^{k-1} \subset \R^k)
= \frac{k\pi^{k/2}}{\Gamma(k/2)}.$$
We combine an explicit computation for 
hemispheres with Lemma \ref{half-space-isop}
to see that 
\begin {eqnarray} \label{isop-ineq1}
\int_{\del \widetilde D} x_n^2 dA & \geq & 
\omega_{n-1}^{\frac{1}{n+2}} \frac{\int_0^{\pi/2}
\sin^{n-2}x \cos^2 x dx}{\left [ \int_0^{\pi/2}
\sin^{n-2} x \cos^{4}x dx \right ]^
{\frac{n+1}{n+2}}} \cdot \left [ 
\int_{\widetilde D}x_n^2 dV \right ]^{\frac{n+1}
{n+2}} \\
& = & c(n) \left [ \int_{\widetilde D} x_n^2 dV 
\right ]^{\frac{n+1}{n+2}} \nonumber
\end {eqnarray}
for all bounded domains $\widetilde D \subset 
\R^n_+$. 

Next we open the wedge $\calW$ up to the 
upper half-space $\R^n_+$. To describe this 
map, recall that the first eigenvalue 
of $\Omega$ is $\mu$, with eigenfunction $\psi$, 
and that our harmonic weight function is 
$w(r,\theta) = r^\alpha \psi(\theta)$, with 
$\alpha$ given by (\ref{defn-harm-weight}). 
Now let
\begin{equation} \label {opening-eqn}
\Psi : \calW \rightarrow \R^n_+, \qquad 
\Psi(r,\theta) = r^{\frac{2\alpha +n -1}
{n+1}}
 \left ( \Pi (\nabla \psi) , \psi \right )
= (x_1, \dots, x_n).
\end {equation}
Here we consider $\nabla \psi$ as a vector 
in $T_\theta \Ss^{n-1} \subset \R^n$, and let 
$\Pi$ be orthogonal projection onto the 
horizontal hyperplane $\{ x_n = 0 \}$. Note 
that, because $\Omega$ is contained in 
the open upper hemisphere, $\Pi:T_\theta
(\Ss^{n-1}) \rightarrow \R^{n-1} = 
\{ x_n = 0\}$ is a 
linear isomorphism for each $\theta \in 
\Omega$, and so $\Pi(\nabla 
\psi) = 0$ only if $\nabla \psi =0$. Also, 
$\Omega$ is convex, so a theorem of Brascamp 
and Lieb \cite{BL} (also see \cite{Kor}) implies 
the first  eigenfunction $\psi$ is log concave. 
In particular, the only critical point of 
$\psi$ is its maximum, and this is a nondegenerate 
critical point. 
We conclude $\Psi$ is indeed a diffeomorphism. 

Next we estimate the effect of $\Psi$ on 
the arclenth and area elements. Let $d\theta$ 
be the standard arclength element on $\Ss^{n-1}$, 
so that $ds = dr^2 + r^2 d\theta^2$ is the 
arclength element in $D$. We denote the arclenth 
element in the image domain $\widetilde D$ by 
$d\widetilde s^2 = dx_1^2 + \cdots + dx_n^2$. 
We also denote the corresponding area elements by 
$dA$ (on $\del D$) and $d\widetilde A$ (on $\del 
\widetilde D$). A quick computation gives 
\begin {eqnarray*}
dx_1^2 + \cdots + dx_n^2 & \leq & 
\left ( \frac{2\alpha +n-1}{n+1} 
\right )^2 r^{\frac{4\alpha -4}{n+1}}
(\psi^2 + |\nabla \psi|^2) dr^2 + 
r^{\frac{4\alpha +2n -2}{n+1}}
(|\nabla \psi|^2 + |D^2 \psi|^2) d\theta^2 \\
& \leq & c \| \psi\|^2_{C^2} r^{\frac{4\alpha 
-4}{n+1}} (dr^2 + r^2 d\theta^2 ).
\end {eqnarray*}
Here $|D^2 \psi|$ is the $L^2$-norm of $D^2\psi$, 
considered as a quadratic form. We have also 
used the fact that $\Pi$ is norm-decreasing, so 
(for instance) $|\Pi(\nabla \psi)| \leq |\nabla \psi|$. 
Combine the computation above with estimates 
for the eigenfunction $\psi$ (see \cite{Li}) 
to conclude that there is a constant $\hat c 
= \hat c(n, \Omega)$ such that
$$dx_1^2 + \cdots dx_n^2 \leq \hat c(n,\Omega)
r^{\frac{4(\alpha-1)}{n+1}} (dr^2 + r^2 d\theta^2).$$
Therefore, for some constant $\hat c$ depending 
on $n$ and $\Omega$, 
\begin {eqnarray} \label {area-comp}
x_n^2 d\widetilde A & \leq & \hat c
r^{\frac{4\alpha + 2n -2}{n+1}} \psi^2
\left ( r^{\frac{2(\alpha -1)}{n+1}}
\right )^{n-1} dA 
= \hat c r^{2\alpha} \psi^2 dA
\\ \nonumber 
& = & \hat c w^2 dA.
\end {eqnarray}

In order to prove the inequality 
(\ref{weight-isop}) we need one final 
ingredient, which is an inequality of 
Szeg\H{o} \cite{Sz} (also see Lemma 2 
of \cite{RT}). 
If $f\geq 0$, $F' = f$, $g$ is nondecreasing, 
and $G'= g$, then for any bounded measurable 
set $E\subset \R$,
\begin {equation} \label {szego-ineq}
G \left ( \int_E f(x) dx \right ) 
\leq \int_E g(F(x)) f(x) dx, \end{equation}
with equality if and only if $E$ is almost 
everywhere an interval of the form $[0,R]$. 

We now prove our weighted isoperimetric 
inequality.  
\begin{proof}
Let $D$ be a bounded domain in the wedge 
$\calW$, and $\widetilde D = \Psi(D)$ its 
image under the change of variables. We 
apply inequality (\ref{area-comp}) 
to the isoperimetric inequality 
(\ref{isop-ineq1}), and conclude
\begin {eqnarray} \label {isop-ineq3/2}
\int_{\del D} w^2 dA & \geq & c 
\int_{\del \widetilde D} x_n^2 d\widetilde A 
\geq	c 
\left [ \int_{\widetilde D} x_n^2 
d\widetilde V \right 
]^{\frac{n+1}{n+2}} \\ \nonumber 
& = & c \left [ \int_D r^{\frac{4\alpha +2n-2}
{n+1}} \psi^2 |\det(D\Psi)| dr dA_{\Ss^{n-1}} 
\right ]^{\frac {n+1}{n+2}} \\ \nonumber 
& \geq & c \left [ \int_D r^{\frac{4\alpha 
+ 2n -2}{n+1}} \psi^2 r^{\frac{2\alpha - 2}
{n+1}} \left ( r^{\frac{2\alpha +n-1}{n+1}}
\right )^{n-1} dr dA_{\Ss^{n-1}} 
\right ]^{\frac{n+1}{n+2}} \\ \nonumber
& = & c \left [ \int_D  r^{\frac{2n\alpha 
+ n^2 + 4\alpha  -3}{n+1}} \psi^2 dr 
dA_{\Ss^{n-1}}\right ]^{\frac{n+1}{n+2}}.
\end {eqnarray}
Here $c$ is a generic constant depending only 
on $n$ and $\Omega$ which we allow to change 
line to line in the computation. The estimate 
of $|\det D\Psi|$ arises from 
$$D\Psi = \left [ \begin {array}{cc} 
r^{\frac{2\alpha + n -1}{n+1}} D(\Pi(\nabla 
\psi)) & \left ( \frac{2\alpha + n -1}
{n+1}\right ) r^{\frac{ 2\alpha -2}{n+1}} \Pi 
(\nabla \psi) \\
r^{\frac{2\alpha +n -1}{n+1}} \Pi 
(\nabla \psi) & \left ( \frac{2\alpha 
+ n-1}{n+1} \right ) r^{\frac{2\alpha 
-2}{n+1}} \psi \end {array} \right ]. $$
Restricted to $\{ r= 1\}$, the map $\Psi$ is 
a diffeomorphism between compact 
domains, and so there are positive 
constants $c_1 < c_2$ such that 
$$c_1 \leq |\left . \det D\Psi\right |_{r=1}|
\leq c_2.$$
Now compute $\det D\Psi$ in general by 
expanding along the last column. Each entry 
in the last column of $D\Psi$ has a factor 
or $r^{\frac{2\alpha -2}{n+1}}$, while all the 
other entries have a factor of $r^{\frac
{2\alpha +n-1}{n+1}}$. Thus 
$$|\det D\Psi| \geq c_1 r^{\frac{2\alpha 
-2}{n+1}} \left ( r^{\frac{2\alpha +n -1}{n+1}}
\right )^{n-1} = c_1 r^{\frac{2n\alpha + n^2 
-2n-1}{n+1}}.$$ 

Denote the inner integral of 
(\ref{isop-ineq3/2}) by 
$$I = \int_D r^{\frac{2n\alpha +n^2 +
4\alpha -3}{n+1}} \psi^2 dr dA_{\Ss^{n-1}}.$$
For a given $\theta \in \Omega$, we let 
$L_\theta$ be the radial slice 
$$L_\theta = \{ r \in [0,\infty) 
\mbox{ }|\mbox{ } (r,\theta) \in D \}, $$
and rewrite the integral $I$ as 
$$I = \int_\Omega \int_{L_\theta}  
r^{\frac{2n\alpha + n^2 
+ 4\alpha -3}{n+1}} \psi^2 dr dA_{\Ss^{n-1}} 
= \int_\Omega \int_{L_\theta} r^{\beta}
(r^{\beta+ 1})^\gamma dr \psi^2(\theta)
dA_{\Ss^{n-1}}. $$ 
Here $\beta$ and $\gamma$ are positive 
parameters we will choose later. 

We let  
$$f(r) = (\beta +1) r^\beta, \qquad 
g(x) = x^\gamma,$$
and apply inequality (\ref{szego-ineq}), 
yielding 
\begin {eqnarray} \label {isop-ineq2}
I & = & \int_\Omega \int_{L_\theta} 
r^\beta (r^{\beta+1})^\gamma dr \psi^2
(\theta) dA_{\Ss^{n-1}} \\ \nonumber 
& = & \frac{1}{\beta+1}
\int_\Omega \int_{L_\theta} f(r) g(F(r)) dr 
\psi^2 dA_{\Ss^{n-1}} \\ \nonumber
& \geq & \frac{1}{\beta+1} \int_\Omega 
G \left ( \int_{L_\theta} f(r) dr
\right ) \psi^2(\theta) dA_{\Ss^{n-1}} 
\\ \nonumber
& = & \frac{1}{(\beta+1)(\gamma+1)} 
\int_{\Omega} \left ( \int_{L_\theta}
r^\beta dr \right )^{\gamma+1} \psi^2(\theta)
dA_{\Ss^{n-1}}.
\end {eqnarray}
Moreover, equality only occurs if $L_\theta$ is 
an interval of the form $(0,R(\theta))$ for 
almost every $\theta$. 

Define the measure $\nu$ on $\Omega$ by 
$d\nu = \psi^2 dA_{\Ss^{n-1}}$, and notice 
it has  total measure $\nu(\Omega) =1$. 
Applying H\"older's inequality with exponents 
$\gamma+1$ and $\frac{\gamma+1}{\gamma}$, 
to the functions $h(\theta) = \int_{L_\theta}
r^\beta dr$ and $1$, gives us 
\begin {eqnarray*}
\left [ \int_\Omega \left ( \int_{L_\theta}
r^\beta dr \right )^{\gamma+1} d\nu 
\right ]^{\frac{1}{\gamma+1}} & = & 
 \left [ \int_\Omega \left ( \int_{L_\theta}
r^\beta dr \right )^{\gamma+1} d\nu 
\right ]^{\frac{1}{\gamma+1}} \left [
\int_\Omega d\nu \right ]^{\frac{\gamma}
{\gamma + 1}} \\
& \geq & \int_\Omega \int_{L_\theta} r^\beta
dr d\nu.
\end {eqnarray*}
Raise each side of this inequality 
to the power $\gamma+1$, and plug 
the result into inequality (\ref{isop-ineq2}) 
to obtain 
\begin {equation} \label{isop-ineq3}
I \geq \frac{1}{(\beta + 1)(\gamma + 1)}
\left ( \int_\Omega \int_{L_\theta}
r^\beta dr \psi^2(\theta) dA_{\Ss^{n-1}}
\right )^{\gamma+1}. \end {equation}
Moreover, equality only occurs if 
$\int_{L_\theta} r^\beta dr$ is constant 
in $\theta$. 

Plug the inequality (\ref{isop-ineq3}) into 
(\ref{isop-ineq3/2}) to get
\begin {equation} \label{isop-ineq4} 
\int_{\del D} w^2 dA \geq 
c\left ( \int_\Omega \int_{L_\theta}
r^\beta dr \psi^2 dA_{\Ss^{n-1}} \right 
)^{\frac{(\gamma+1)(n+1)}{n+2}}.
\end {equation}
Finally, we choose the parameters $\beta$ and 
$\gamma$. In order for the integral to be 
$\int_D w^2 dV$, we need $\beta = 2\alpha + n
- 1$. Thus
$$r^{\frac{2n\alpha + 2n^2 + 4\alpha -3}{n+1}}
= r^\beta (r^{\beta+1})^\gamma \Rightarrow 
\gamma = \frac{2\alpha -2}{2n\alpha + n^2 + 
2\alpha + n} = \frac{2\alpha -2}{(2\alpha + n)
(n+1)}.$$
Plugging these choices of $\beta$ and $\gamma$ 
into inequality (\ref{isop-ineq4}) gives 
us 
$$\int_{\del D} w^2 dA \geq c \left [ 
\int_D w^2 dV \right ]^{\frac{2\alpha + n -1}
{2\alpha + n}} = c \left [ \int_D
w^2 dV \right ]^{\frac{2a+1}{2a+2}},$$
as claimed. 

Now consider the case of equality. To achieve 
equality, we must have equality in 
(\ref{isop-ineq2}), which forces $L_\theta 
= \{r \in [0,\infty) \mbox{ }|\mbox{ } 
(r,\theta) \in D \} = (0,R(\theta))$ almost 
everywhere. Moreover, equality 
in our use of H\"older's inequality forces 
$\int_{L_\theta} r^\beta dr$ to be constant 
in $\theta$. If $L_\theta = (0,R(\theta))$, 
this latter integral 
is $\frac{1}{\beta+1} R^{\beta+1}(\theta)$, 
which can only be constant if the function 
$R(\theta)$ is constant. The combination of 
these two properties forces $D$ to be a 
sector.  

At this point, we can recover the constant $c$
by examining the case of the sector. If $D = 
\calS_{r_0}$, then 
$$\int_{\calS_{r_0}} w^2 dV  = 
\frac{1}{2a+2} r_0^{2a+2}, 
\qquad \int_{\del \calS_{r_0}} w^2 dA = 
r_0^{2a+1}.$$ 
Solving for $r_0$ and rearranging, we get 
$$\int_{\del \calS_{r_0}} w^2 dA = 
\left [ (2a+2) \int_{\calS_{r_0}} w^2 
dV \right ]^{\frac{2a+1}{2a+2}},$$
and so 
$$c = (2a+2)^{\frac{2a+1}{2a+2}}$$
as we claimed. 
\end {proof}

\section{Estimating the Rayleigh quotient}
\label {rayleigh-sec}

In this section we prove our main theorem 
by estimating the Rayleigh quotient 
$$\frac{\int_D |\nabla u|^2 dV}
{\int_D u^2 dV}$$
for an appropriate test function. We proceed 
as in \cite{PW} and \cite{RT}.

\begin {proof}
Let $u \geq 0$, and write $u(r,\theta) = 
w(r,\theta) v(r,\theta)$ for some $v \in 
C^2_0 (\calW)$. Using the fact that $\Delta w 
= 0$ and the divergence theorem, we see 
\begin{eqnarray*}
\int_D vw \langle \nabla v, \nabla w \rangle dV
& = & - \int_D v\Div (vw \nabla w)  dV
= - \int_D v(vw \Delta w + v |\nabla w|^2 
 + w \langle \nabla v, \nabla w\rangle dV \\
 & = & -\int_D v^2 |\nabla w|^2 + vw \langle 
 \nabla v, \nabla w \rangle dV.
\end {eqnarray*}
Rearranging yields 
$$2\int_D vw \langle \nabla v, \nabla w \rangle dV
= - \int_D v^2 |\nabla w|^2 dV,$$
and so 
$$\int_D |\nabla u|^2 dV = \int_D 
v^2 |\nabla w|^2 + 2vw \langle \nabla v, \nabla w
\rangle + w^2 |\nabla v|^2 dv = \int_D w^2 
|\nabla v|^2 dV.$$

For $0 \leq t \leq \bar v = \max(v)$, let 
$D_t := v^{-1}((t,\bar v])$ and 
$$\zeta(t) := \int_{D_t} w^2 dV.$$
By the coarea formula, 
$$\frac{d\zeta}{dt} = -\int_{\del D_t} 
\frac{w^2}{|\nabla v|} dA < 0,$$
and so $\zeta$ is a decreasing function 
and has an inverse $t = t(\zeta)$. Then 
\begin {equation} \label{coarea-1}
\int_D w^2 v^2 dV = \int_0^{\bar v}
t^2 \int_{\del D_t} \frac{w^2}{|\nabla 
v|} dA dt = - \int_0^{\bar v} t^2 
\frac{d\zeta}{dt}(t)dt 
= \int_0^{\bar \zeta} t^2 d\zeta.
\end{equation}
Here $\bar \zeta = \zeta(0) = \max(\zeta)$.

Next we estimate the Dirichlet energy. For this 
computation, it will be convenient to define 
$$p = \frac{2a+1}{2a+2} = \frac{2\alpha + n-1}
{2\alpha +n}.$$ 
Note that, applying the Cauchy-Schwartz 
inequality to $w^2 = (w |\nabla v|^{-1/2})
(w |\nabla v|^{1/2})$, we have 
$$\int_{\del D_t} w^2 dA \leq \left ( 
\int_{\del D_t} \frac{w^2}{|\nabla v|} dA 
\right )^{1/2} \left ( \int_{\del D_t} 
w^2 |\nabla v| dA \right )^{1/2} 
\quad \Rightarrow \quad 
\int_{\del D_t} w^2 |\nabla v| dA
\geq \frac{\left ( \int_{\del D_t} 
w^2 dA \right )^2}{ \int_{\del D_t}
\frac{w^2}{|\nabla v|} dA}.$$
Now use the coarea formula, the 
Cauchy-Schwartz inequality as above, and 
our inequality (\ref{weight-isop}) to 
conclude
\begin{eqnarray} \label{rayleigh-q1}
\int_D |\nabla u|^2 dV & = & 
\int_D w^2 |\nabla v|^2 dV = 
\int_0^{\bar v} \int_{\del D_t} w^2 |\nabla 
v| dA dt \\ \nonumber
& \geq & \int_0^{\bar v} \frac { \left (
\int_{\del D_t} w^2 dA \right )^2}
{\int_{\del D_t} \frac{w^2}{|\nabla v|}dA} dt 
\\ \nonumber 
& \geq & (2a+2)^{\frac{2a+1}{a+1}} 
\int_0^{\bar v} \frac {\left (
\int_{D_t} w^2 dV \right )^{2p}}
{-\frac{d\zeta}{dt}} dt \\ \nonumber 
& = & (2a+2)^{\frac
{2a+1}{a+1}} \int_0^{\bar v}
\frac{\zeta^{2p}(t)}{-\frac{d\zeta}{dt}} dt 
\\ \nonumber
& = & (2a+2)^{\frac{2a+1}{a+1}} \int_0^{\bar \zeta} 
\zeta^{2p} \left ( 
\frac{dt}{d\zeta} \right )^2 d\zeta .
\end {eqnarray} 

Now we define $\lambda_*$ to be the least 
eigenvalue of
\begin {equation} \label{comp-ode} 
\frac{d}{d\zeta} \left ( \zeta^{2p}
\frac{dt}{d\zeta} \right ) + \lambda_* t = 0, 
\qquad t(\bar \zeta) = 0 = \lim_{\zeta 
\rightarrow 0^+} \left (\zeta^{2p} 
\frac{dt}{d\zeta} \right ). \end {equation}
Observe that, by the boundary conditions 
in equation (\ref{comp-ode}), 
$$\int_0^{\bar \zeta} \zeta^{2p} 
\left ( \frac{dt}{d\zeta} \right )^2 d\zeta
= -\int_0^{\bar \zeta} t \frac{d}{d\zeta}
\left ( \zeta^{2p} \frac{dt}{d\zeta} \right ) 
d\zeta.$$
Thus we combine inequalities (\ref{coarea-1}) 
and (\ref{rayleigh-q1}) with the eigenvalue 
equation (\ref{comp-ode}) to get 
\begin {equation} \label{rayleigh-q2}
(2a+2)^{\frac{2a+1}{a+1}} \lambda_* 
= (2a+2)^{\frac{2a+1}{a+1}} \inf \left [ 
\frac{\int_0^{\bar \zeta} \zeta^{2p} \left (
\frac{dt}{d\zeta} \right )^2 d\zeta}
{\int_0^{\bar \zeta} t^2 d\zeta} \right ]
\leq \inf \left [ \frac{\int_D |\nabla u|^2
dV}{\int_D u^2 dV} \right ] = \lambda_1(D).
\end {equation}

We can write the solution to (\ref{comp-ode}) in 
terms of Bessel functions. We try 
$$t(\zeta) = \zeta^q J_b \left ( \frac{\sqrt
{\lambda_*}}{m} \zeta^m \right ),$$
where $J_b$ is the Bessel function of weight 
$b$. Taking a derivative, we see 
$$t' = q \zeta^{q-1} J_b + \sqrt{\lambda_*} 
\zeta^{q+m-1} J_b',$$
and so equation (\ref{comp-ode}) becomes 
\begin {eqnarray*}
-\lambda_* \zeta^q J_b &= &\frac{d}{d\zeta} 
(\zeta^{2p} t'(\zeta)) = 
\frac{d}{d\zeta} (q\zeta^{2p+q-1} J_b + 
\sqrt{\lambda_*}\zeta^{2p+q+m-1} J_b') \\
& = & (2pq + q^2 -q)\zeta^{2p+q-2} J_b 
+ \sqrt{\lambda_*}(2p+2q+m-1) \zeta^{2p+q+m-2} 
J_b' + \lambda_* \zeta^{2p+q+2m-2} J_b'' \\ 
& = & (b^2m^2 +2pq+q^2-q)\zeta^{2p+q-2} J_b - \lambda_* 
\zeta^{2p+q+2m-2} J_b + \sqrt{\lambda_*} (2p+2q-1)
\zeta^{2p+q+m-2} J_b',
\end {eqnarray*}
where we have used Bessel's identity in the last 
equality. We force the last term to be zero by choosing 
$$q = \frac{1-2p}{2} = \frac{2-2\alpha-n}
{4\alpha +n} = -\frac{2a}{4a-n+4}.$$
We force the middle term to be equal to the left 
hand side $-\lambda_* t = -\lambda_* \zeta^q J_b$ 
by choosing 
$$m=1-p = \frac{1}{2\alpha + n} = 
\frac{1}{2a+2}.$$
Finally, we force the last term to be zero by choosing 
$$b^2 = \left ( \frac{2p-1}{2p-2} \right )^2.$$
Notice that we have a choice in the sign 
of the square root, and we choose 
\begin {eqnarray*}
b & = & -\frac{2p-1}{2p-2} = \frac{1-2p}{2p-2}\\
& = & \frac{ 1 - \frac{2a+1}{a+1}} 
{\frac{2a+1}{2a+1} - 2} = 
\frac{a+1-2a-1}{2a+1-2a-2} \\
& = & a.
\end {eqnarray*}
Putting this all together, we have 
$$t(\zeta) = \zeta^{-\frac{2a}{4a-n+4}}
J_a \left ( (2a+2)\sqrt{\lambda_*} 
\zeta^{\frac{1}{2a+2}} \right ).$$

Next we use the boundary condition $t(\bar \zeta) 
= 0$ to find $\lambda_*$. If $j_a$ is the first positive 
zero of $J_a$, we must have 
$$j_a = (2a+2) \sqrt{\lambda_*} \bar 
\zeta^{\frac{1}{2a+2}}$$
and so 
$$
\lambda_* = \frac{\bar \zeta^{-\frac{1}{a+1}}}
{(2a+2)^2} j_a^2 = \frac{1}{(2a+2)^2}
\left [ \int_D w^2 dV \right ]^{-\frac{1}{a+1}}
j_a^2.$$
Using the explicit computation of $\lambda_*$ 
above to compare inequality (\ref{rayleigh-q2}) 
with Lemma \ref{eigen-sector} completes 
the proof of our main theorem. 
\end {proof}

\end{document}